\documentclass[11pt,leqno]{article}
\usepackage{amssymb}
\usepackage[mathscr]{eucal}
\usepackage{amsmath,amssymb,latexsym,theorem,bbm,graphics,pdfpages}
\usepackage{color}
\usepackage{appendix}
\usepackage{epsfig,rotating}

\newcommand{\NN}{\mathbb{N}}

\newcommand{\RR}{\mathbb{R}}

\newcommand{\ZZ}{\mathbb{Z}}

\newcommand{\cB}{{\mathcal B}}

\newcommand{\cF}{{\mathcal F}}

\newcommand{\dd}{\mathrm{d}}

\newcommand{\EE}{\operatorname{\mathbb{E}}}
\newcommand{\PP}{\operatorname{\mathbb{P}}}

\renewcommand{\mid}{\,|\,}

\renewcommand{\leq}{\leqslant}
\renewcommand{\geq}{\geqslant}

\newcommand{\proofend}{\hfill\mbox{$\Box$}}

\numberwithin{equation}{section}

\theoremstyle{change} \theorembodyfont{\em}
\newtheorem{Lem}{Lemma.}[section]
\newtheorem{Thm}{Theorem.}[section]
\newtheorem{Pro}{Proposition.}[section]
\newtheorem{Cor}{Corollary.}[section]
\newtheorem{Def}{Definition.}[section]

\theorembodyfont{\rm}
\newtheorem{Rem}{Remark.}[section]
\newtheorem{Ex}{Example.}[section]

\begin{document}

\begin{center}
 {\bfseries\Large A Robbins-Monro type algorithm \\ for computing global minimizer of generalized conic functions} \\[5mm]

 {\sc\large M\'aty\'as $\text{Barczy}^{*,\diamond}$,
            \ \'Abris $\text{Nagy}^{**}$, \ Csaba $\text{Nosz\'aly}^{*}$, \ Csaba $\text{Vincze}^{***}$}
\end{center}

\vskip0.2cm

\noindent * Faculty of Informatics, University of Debrecen,
            P.~O.~Box 12, 4010 Debrecen, Hungary.
            Tel.:  +36-52-512900, Fax: +36-52-512996

\noindent ** Institute of Mathematics, MTA-DE Research Group ''Equations Functions and Curves'',
             Hungarian Academy of Sciences and University of Debrecen,
             P.~O.~Box 12, 4010 Debrecen, Hungary.

\noindent *** Institute of Mathematics, University of Debrecen,
              P.~O.~Box 12, 4010 Debrecen, Hungary.

\noindent e--mails: barczy.matyas@inf.unideb.hu (M. Barczy),
                    abris.nagy@science.unideb.hu (\'A. Nagy),\\
                    noszaly.csaba@inf.unideb.hu (Cs. Nosz\'aly),
                    csvincze@science.unideb.hu (Cs. Vincze).

\noindent $\diamond$ Corresponding author.


\vspace*{0.5cm}

\renewcommand{\thefootnote}{}
\vspace*{0.2cm}
\footnote{
The research of M\'aty\'as Barczy was realized in the frames of
 T\'AMOP 4.2.4.\ A/2-11-1-2012-0001 ,,National Excellence Program --
 Elaborating and operating an inland student and researcher personal support
 system''.
The project was subsidized by the European Union and co-financed by the
 European Social Fund.
\'Abris Nagy has been supported, in part, by the Hungarian Academy of
 Sciences, the European Union and the State of Hungary, co-financed by
 the European Social Fund in the framework of T\'AMOP 4.2.4. A/2-11-1-2012-0001 `National Excellence Program'.
Csaba Vincze was partially supported by the European Union and the
 European Social Fund through the project Supercomputer, the national
 virtual lab (grant no.: T\'AMOP-4.2.2.C-11/1/KONV-2012-0010).
Csaba Vincze is supported by the University of Debrecen's internal
research project RH/885/2013.}

\vspace*{-15mm}

\begin{abstract}
We generalize the notion and some properties of the conic function introduced by Vincze and Nagy (2012).
We provide a stochastic algorithm for computing the global minimizer of generalized conic functions,
 we prove almost sure and \ $L^q$-convergence of this algorithm.
\end{abstract}

\section{Introduction}

Let \ $K$ \ be a compact body in \ $\RR^2$ \ (a non-empty compact set coinciding with the closure of its interior)
 and consider the distance function induced by the taxicab norm.
The so called conic function \ $F_K$ \ associated to \ $K$ \ (introduced by Vincze and Nagy \cite[Defintion 6]{VinNag}, see also
 Definition \ref{Def_gen_conic_VN}) measures the average taxicab distance of the points from \ $K$ \ via integration
 with respect to the Lebesgue measure, or explaining in another way: the conic function \ $F_K$ \ at some point
 \ $(x,y)\in\RR^2$ \ can be interpreted as the expectation of the random variable defined as the taxicab distance of
 \ $(x,y)$ \ and \ $(\xi,\eta)$, \ where \ $(\xi,\eta)$ \ is a uniformly distributed random variable on \ $K$, \ for more details
 see part (ii) of Remark \ref{Rem1}.
Conic functions are extensively used in geometric tomography since they contain a lot of information about unknown bodies, for a more detailed
 discussion see Gardner \cite{Gar} and Vincze and Nagy \cite{VinNag}.
We call the attention that in the literature one can find other definitions of ''conic functions''
 that are completely different from ours.
For example, in optimization a conic function is usually defined to be the ratio of a quadratic function and the square of a linear function
 on the open halfspace, where the linear function is positive, see, e.g., Luksan \cite[formula (2.1)]{Luk}.
Wang et al. \cite{WanWanRen} introduced another definition of conic functions in metric spaces and obtained a new condition
 for metric spaces being compact in terms of conic functions.

We recall that one of the striking features of the conic function \ $F_K$ \ is that a point in \ $\RR^2$ \ is a global minimizer
 of \ $F_K$ \ if and only if it bisects the area of \ $K$, \ i.e., the vertical and horizontal lines through this point
 cut the compact body \ $K$ \ into two parts with equal areas, see Vincze and Nagy \cite[Corollary 1]{VinNag}.
We call the attention that points with similar properties are important and well-studied in geometry.
For instance, we mention that if \ $S$ \ is a convex set in \ $\RR^2$, \ then there exist two perpendicular lines that divide
 \ $S$ \ into four parts with equal areas, see Yaglom and Boltyanskii \cite[Section 3]{YagBol}.

In Section \ref{Section_conic} of the present paper we generalize the conic function \ $F_K$ \ introduced by Vincze and Nagy \cite{VinNag}
 in a way that it measures the average taxicab distance of the points from \ $K$ \ via integration with respect to some measure
  \ $\mu$ \ on \ $K$ \ with \ $\mu(K)<\infty$, \ see Definition \ref{Def_gen_conic}.
From geometric point of view the body \ $K$ \ associated with some measure \ $\mu$ \ can be considered as a mathematical model
 of a non-homogeneous body and hence our generalization of conic functions may find applications in (geometric) tomography where
 typically non-homogeneous bodies occur.
We generalize Theorems 3, 4, 5, Lemmas 6, 7 and Corollary 1 in Vincze and Nagy \cite{VinNag} for conic functions \ $F_{K,\mu}$ \
 associated with a compact body \ $K$ \ and a measure \ $\mu$ \ with \ $\mu(K)<\infty$.
\ We only mention that it turns out that a point in \ $\RR^2$ \ is a global minimizer of \ $F_{K,\mu}$ \ if and only if it bisects
 the \ $\mu$-area of \ $K$, \ see Corollary \ref{Cor1}.

In Section \ref{Section_stoch_alg_gen_conic} we give a stochastic algorithm for the global minimizer of the convex function \ $F_{K,\mu}$.
\ In the heart of our algorithm the well-known Robbins-Monro algorithm (see \cite{RobMon}) lies, and we prove almost sure and \ $L^q$-convergence
 of our algorithm.
More precisely, we define recursively a sequence \ $(X_k)_{k\in\ZZ_+}$ \ of random variables
 (see \eqref{recursion}) which forms an inhomogeneous Markov chain
 and we prove almost sure and \ $L^q$-convergence of this Markov chain via Robbins-Monro algorithm,
 see Theorem \ref{Thm3}.
We also prove almost sure and \ $L^q$-convergence of the sequence \ $(F_{K,\mu}(X_k))_{k\in\NN}$, \ see Theorem \ref{Thm4}.
\ In general, stochastic algorithms for finding a minimum of a convex function have a vast literature,
 see, e.g., Robert and Casella \cite{RobCas} and Bouleau and L\'epingle \cite{BouLep}.
Without giving an introduction of the newest results in the field we only mention the paper \cite{ArnDomPhaYan} of  Arnaudon et al.,
 which in some sense motivated our study.
They gave a stochastic algorithm which converges almost surely and in \ $L^2$ \ to the so-called \ $p$-mean of a probability
 measure supported by a regular geodesic ball in a manifold.

\section{Generalized conic functions}\label{Section_conic}

Let \ $\ZZ_+$, \ $\NN$, \ $\RR$ \ and \ $\RR_+$ \ denote the set
 of non-negative integers, positive integers, real numbers and non-negative real numbers, respectively.
For an \ $x\in\RR^2$, \ we will denote its Euclidean norm by \ $\Vert x\Vert$.
\ Let \ $K\subset\RR^2$ \ be a non-empty compact set such that it coincides with the closure of its interior.
In geometry \ $K$ \ is called a compact body.
By \ $\cB(\RR^d)$ \ and \ $\cB(K)$, \ we denote the Borel \ $\sigma$-algebra on \ $\RR^d$ \ and on \ $K$, \ respectively,
 where \ $d\in\NN$.
\ For all \ $x,y\in\RR$ \ let us introduce the following notations
 \begin{align*}
   &\{ K<_1 x\} := \{ (\alpha,\beta)\in K : \alpha < x\},
   \quad\quad\{ x<_1 K\} := \{ (\alpha,\beta)\in K : x < \alpha \},\\
   &\{ K<_2 y\} := \{ (\alpha,\beta)\in K : \beta < y\},
    \quad\quad \{ y<_2 K\} := \{ (\alpha,\beta)\in K : y < \beta \},\\
   &\{ K =_1 x\} := \{ (\alpha,\beta)\in K : \alpha = x\},
   \quad\quad \{ K =_2 y \} := \{ (\alpha,\beta)\in K : \beta = y \}.
 \end{align*}
The notations \ $\{ K\leq_1 x\}$, \ $\{ x\leq_1 K\}$, \ $\{ K\leq_2 y\}$ \ and \ $\{ y\leq_2 K\}$ \ are defined in the same way.
For a function \ $f:\RR^2\to\RR$, \ we will denote by \ $D_1f$ \ and \ $D_2f$ \ the partial derivatives
 of \ $f$.

Next we recall the notion of a generalized conic function associated with \ $K$ \ due to Vincze and Nagy \cite{VinNag}.

\begin{Def}\label{Def_gen_conic_VN}{\bf (Vincze and Nagy \cite[Definition 6]{VinNag})}
The generalized conic function $F_K:\RR^2\to\RR$ associated to $K$ is defined by
 \[
   F_K(x,y):=\frac{1}{A(K)}\int_K d_1((x,y),(\alpha,\beta))\,\dd\alpha\dd\beta,
     \qquad (x,y)\in\RR^2,
 \]
 where \ $A(K)$ \ is the 2-dimensional Lebesgue measure (area) of \ $K$, \ and the distance function
 \ $d_1$ \ is given by \ $d_1((x,y),(\alpha,\beta)):=\vert x-\alpha\vert + \vert y-\beta\vert$,
 $(x,y),(\alpha,\beta)\in\RR^2$ \ ($d_1$ \ is known to be the metric induced by the taxicab norm).
\end{Def}

The next result is about the global minimizer of \ $F_K$.

\begin{Pro}\label{Pro2_VN}{\bf (Vincze and Nagy \cite[Corollary 1]{VinNag})}
A point in \ $\RR^2$ \ is a global minimizer of the generalized conic function \ $F_K$ \ if and only
 if it bisects the area of \ $K$, \ i.e., the vertical and the horizontal lines through this point
 cut the compact body \ $K$ \ into two parts with equal area.
\end{Pro}

We note that the global minimizer of the generalized conic function \ $F_K$ \ is not unique in general.
In Proposition \ref{Pro1} we give a sufficient condition for its uniqueness.

In what follows we will frequently use the following conditions
 \begin{align*}
   &\mathbf{ (C.1)} \qquad\qquad \text{$K$ \ is connected}, \\
   &\mathbf{ (C.2)}\qquad\qquad  \mu(B(p,\varepsilon)\cap K)>0 \quad \text{for all \ $p\in K$, $\varepsilon>0$ \ and \
          $B(p,\varepsilon)$,}
 \end{align*}
 where \ $\mu$ \ is a measure on the measurable space \ $(K,\cB(K))$ \ and
 \ $B(p,\varepsilon)$ \ denotes the open ball around \ $p$ \ with radius \ $\varepsilon$, \ and
 \[
  \mathbf{ (C.3)} \qquad\qquad \mu(\{K=_1x\}) = \mu(\{K=_2y\})= 0 \qquad \text{for all \ $x,y\in\RR$}.
 \]
We call the attention that Condition (C.3) does not hold for a measure in general.
For example, if \ $\mu$ \ is the distribution of a discrete random variable having values in \ $K$,
 \ then Condition (C.3) does not hold.
However, if \ $\mu$ \ is the \ $2$-dimensional Lebesgue measure on \ $K$, \ then
 Conditions (C.2) and (C.3) hold automatically.

\begin{Pro}\label{Pro1}
If Condition (C.1) holds, then the convex function \ $F_K$ \ has a unique global minimizer \ $(x^*,y^*)\in\RR^2$,
 \ that is, \ $F_K(x,y)>F_K(x^*,y^*)$ \ for \ $(x,y)\ne(x^*,y^*)$, \ $(x,y)\in\RR^2$.
\end{Pro}

\noindent{\bf Proof.}
The existence of a global minimizer of \ $F_K$ \ can be checked as follows.
By Theorem 3 in Vincze and Nagy \cite{VinNag}, \ $F_K$ \ is a finite-valued convex function defined on
 \ $\RR^2$ \ and its level sets are compact subsets of \ $\RR^2$.
\ Hence \ $F_K$ \ is continuous and consequently it reaches its minimum on every compact set.

Now we turn to prove the uniqueness of \ $(x^*,y^*)$.
\ Let us suppose that \ $(x^*,y^*)\in\RR^2$ \ and \ $(\widetilde{x^*},\widetilde{y^*})\in\RR^2$ \ are global minimizers
 of \ $F_K$ \ such that \ $(x^*,y^*)\ne (\widetilde{x^*},\widetilde{y^*})$.
\ Then \ $x^*\ne\widetilde{x^*}$ \ or \ $y^*\ne\widetilde{y^*}$.
\ We may assume that \ $\widetilde{x^*}<x^*$.
\ Then both of the vertical lines \ $\RR^2=_1 x^*$ \ and \ $\RR^2=_1 \widetilde{x^*}$ \ bisect the area of \ $K$.
Note that since Condition (C.3) holds automatically for the \ $2$-dimensional Lebesgue measure,
 the bisection of the area of \ $K$ \ is well-defined.
Let us consider the open half-planes
 \[
   H^*:=\RR^2<_1x^* \qquad \text{and}\qquad \widetilde{H^*}:=\RR^2>_1\widetilde{x^*}.
 \]
Note that \ $(\widetilde{x^*},\widetilde{y^*})\in H^*$ \ and \ $(x^*,y^*)\in \widetilde{H^*}$.
\ We show that \ $K\cap (H^* \cap \widetilde{H^*}) = \emptyset$.
\ On the contrary, let us suppose that there exists \ $p\in\RR^2$ \ such that \ $p\in K\cap (H^* \cap \widetilde{H^*})$.
\ Since \ $K$ \ is a non-empty compact body, there exist
 \[
    0<\varepsilon<\min\{ d_2(p,\RR^2=_1 x^*), d_2(p,\RR^2=_1 \widetilde{x^*}) \}
 \]
 and \ $q\in B(p,\varepsilon)$ \ such that \ $q$ \ is an interior point of \ $K$, \ where \ $d_2$ \ denotes
 the standard Euclidean distance on \ $\RR^2$.
\ Hence there exists
 \[
     0<\delta<\min\{ d_2(p,\RR^2=_1 x^*), d_2(p,\RR^2=_1 \widetilde{x^*}) \}
 \]
 such that \ $B(q,\delta)\subset K\cap (H^* \cap \widetilde{H^*})$.
\ Then
 \begin{align}\label{seged20}
  \begin{split}
   &A(K<_1 \widetilde{x^*}) = A(\widetilde{x^*}<_1 K)
       \geq A(B(q,\delta)) + A(x^* <_1 K),\\
   &A(x^* <_1 K) = A(K <_1 x^*)
       \geq A(B(q,\delta)) + A( K <_1 \widetilde{x^*}),
   \end{split}
 \end{align}
 and hence
 \[
    A(K <_1 x^*)\geq 2A(B(q,\delta)) + A(K <_1 x^*),
 \]
 i.e., \ $0\geq A(B(q,\delta))$, \ which yields us to a contradiction.
At this point we implicitly used that  Condition (C.2) holds automatically for the \ $2$-dimensional Lebesgue measure.
Hence \ $K\cap (H^* \cap \widetilde{H^*})=\emptyset$.
\ Let  \ $0<\eta<(x^* - \widetilde{x^*})/2$, \ and let us consider the open half-planes
 \[
   I^*:=\RR^2>_1x^*-\eta \qquad \text{and}\qquad \widetilde{I^*}:=\RR^2<_1\widetilde{x^*}+\eta.
 \]
Then \ $I^*$ \ and \ $\widetilde{I^*}$ \ are open sets of \ $\RR^2$, \ $I^*\cap \widetilde{I^*} = \emptyset$, \ and, since
 \ $K\cap (H^* \cap \widetilde{H^*})=\emptyset$, \ we have \ $K\subset I^*\cup \widetilde{I^*}$.
\ Further, \ $I^*\cap K$ \ and \ $\widetilde{I^*}\cap K$ \ are separated sets such that their union equals \ $K$.
\ This is a contradiction due to the connectedness of \ $K$.
\ Hence \ $x^*= \widetilde{x^*}$, \ and in a similar way we have \ $y^*= \widetilde{y^*}$.
\proofend

\begin{figure}[h!]
\centering
\includegraphics[height=3cm]{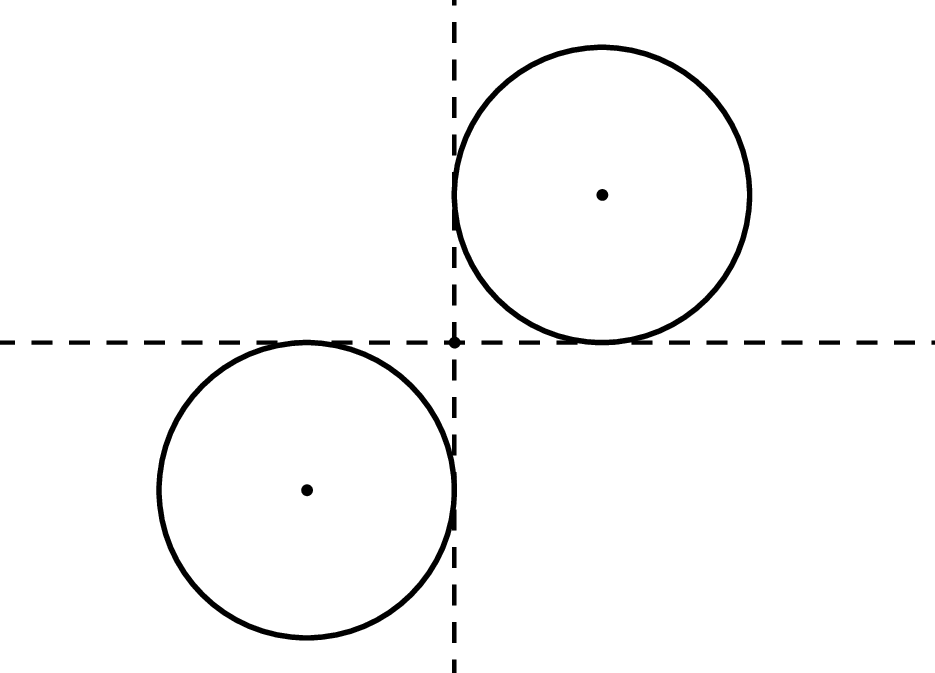}
\hskip1cm
\includegraphics[height=3cm]{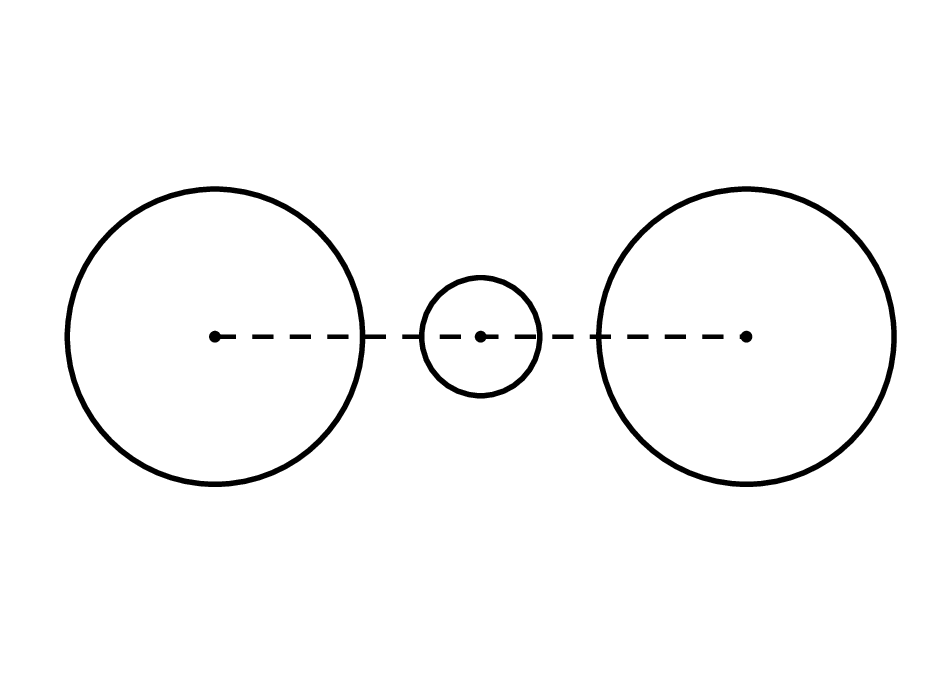}
\hskip1cm
\includegraphics[height=3cm]{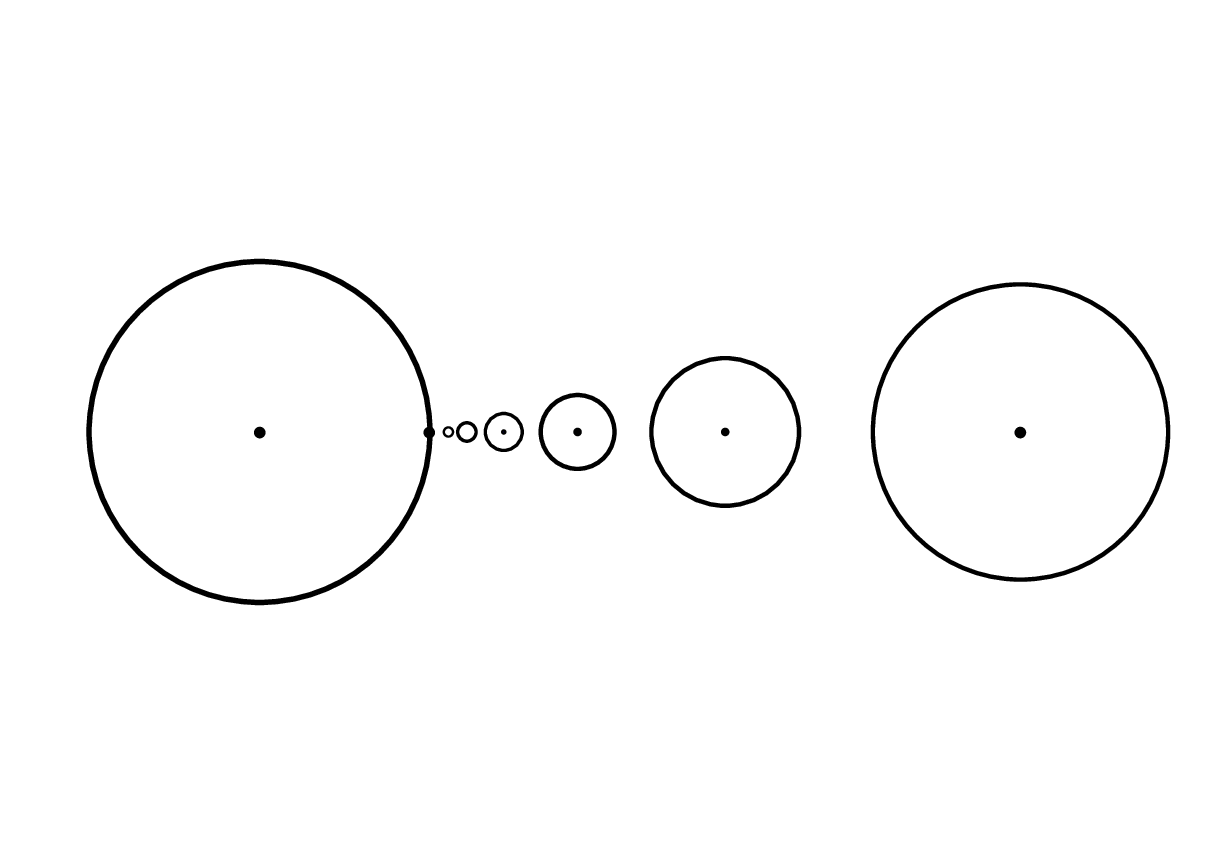}
\caption{Examples for \ $K$ \ such that Condition (C.1) does not hold but \ $F_K$ \ has a unique global minimizer.}
\label{figure1}
\end{figure}

We call the attention that Condition (C.1) is sufficient but not necessary in order that
 the generalized conic function \ $F_K$ \ should have a uniquely determined global minimizer.
 Figure \ref{figure1} shows three different cases where Condition (C.1) is not satisfied but \ $F_K$ \ has a unique
 global minimizer.

On the subfigure (c) of Figure \ref{figure1}, the circles have centers \ $(-1/\sqrt{12},0)$ \
 and \ $(1/2^n,0)$ \ with radii \ $1/\sqrt{12}$ \ and \ $1/2^{n+2}$, \ respectively, where \ $n\in\ZZ_+$.

\begin{Ex}\label{example1}
{(i)} If \ $K$ \ is the square with vertexes \ $(0,0), (0,1), (1,0), (1,1)$, \ then
 \[
    F_K(x,y)=\left(x-\frac{1}{2}\right)^2 + \left(y-\frac{1}{2}\right)^2 + \frac{1}{2},
    \qquad (x,y)\in K,
 \]
 see, e.g., Vincze and Nagy \cite[Example 3]{VinNag}.
Using that \ $K$ \ is connected, by Propositions \ref{Pro2_VN} and \ref{Pro1},
 the global minimizer of \ $F_K$ \ is \ $(x,y)=(\frac{1}{2},\frac{1}{2})$.

\noindent{(ii)} If \ $K$ \ is the triangle with vertexes \ $(0,0), (0,1), (1,0)$, \ then
 \[
    F_K(x,y)=-\frac{2}{3}(x^3+y^3) + 2(x^2+y^2) - (x+y) + \frac{2}{3},
    \qquad (x,y)\in K.
 \]
Indeed, \ $F_K(x,y) = \EE(\vert\xi-x\vert) + \EE(\vert\eta-y\vert)$ \ for all \ $(x,y)\in \RR^2$, \ where \ $(\xi,\eta)$ \ is a
 uniformly distributed random variable on \ $K$.
\ Then the joint density function of \ $(\xi,\eta)$, \ and the density functions of the marginals of \ $(\xi,\eta)$ \ take the forms
 \[
  f_{(\xi,\eta)}(\alpha,\beta)
     =\begin{cases}
        2 & \text{if \ $(\alpha,\beta)\in K$,}\\
        0 & \text{if \ $(\alpha,\beta)\not\in K$,}
      \end{cases}
 \]
 and
 \[
  f_{\xi}(\alpha)
     =\begin{cases}
        -2\alpha+2 & \text{if \ $\alpha\in [0,1]$,}\\
        0 & \text{if \ $\alpha\not\in [0,1]$,}
      \end{cases}
    \qquad \qquad
     f_{\eta}(\beta)
     =\begin{cases}
        -2\beta+2 & \text{if \ $\beta\in [0,1]$,}\\
        0 & \text{if \ $\beta\not\in [0,1]$,}
      \end{cases}
 \]
 respectively.
Hence for all \ $(x,y)\in K$,
 \begin{align*}
  \EE(\vert \xi - x \vert)
    & = \int_0^1 \vert \alpha - x \vert (-2\alpha + 2)\,\dd\alpha\\
    &  = \int_0^x (x-\alpha) (-2\alpha + 2)\,\dd\alpha
        + \int_x^1 (\alpha-x) (-2\alpha + 2)\,\dd\alpha\\
    & = -\frac{2}{3}x^3 + 2x^2 - x + \frac{1}{3},
 \end{align*}
 and similarly \ $\EE(\vert \eta - y \vert) = -\frac{2}{3}y^3 + 2y^2 - y + \frac{1}{3}$ \ for all \ $(x,y)\in K$.
\ Hence the global minimizer of \ $F_K$ \ is \ $(1-\sqrt{2}/2,1-\sqrt{2}/2)$.
\ Indeed, the solution in \ $K$ \ of the system of equations
 \begin{align*}
   D_1F_K(x,y) = -2x^2 + 4x -1=0\qquad\text{and}\qquad
   D_2F_K(x,y) = -2y^2 + 4y -1=0,
 \end{align*}
 is \ $(1-\sqrt{2}/2,1-\sqrt{2}/2)$.
Using that \ $K$ \ is connected, by Propositions \ref{Pro2_VN} and \ref{Pro1},
 the global minimizer of \ $F_K$ \ is \ $(1-\sqrt{2}/2,1-\sqrt{2}/2)$.
\end{Ex}

In what follows we generalize the notion of the conic function introduced
 by Vincze and Nagy \cite[Definition 6]{VinNag}, see also Definition \ref{Def_gen_conic_VN}.

\begin{Def}\label{Def_gen_conic}
Let \ $\mu$ \ be a measure on the measurable space \ $(K,\cB(K))$ \ such that \ $\mu(K)<\infty$.
\ The generalized conic function $F_{K,\mu}:\RR^2\to\RR$ associated to \ $K$ \ and \ $\mu$ \ is defined by
 \[
   F_{K,\mu}(x,y):=\int_K d_1((x,y),(\alpha,\beta))\,\mu(\dd\alpha,\dd\beta),
     \qquad (x,y)\in\RR^2.
 \]
\end{Def}

\begin{Rem}\label{Rem1}
{(i):} Note that under the conditions of Definition \ref{Def_gen_conic} we have \ $F_{K,\mu}(x,y)$ \ is well-defined for all
 \ $(x,y)\in\RR^2$, \ since for fixed \ $(x,y)\in\RR^2$, \ the function \ $K\ni(\alpha,\beta)\mapsto d_1((x,y),(\alpha,\beta))$ \
 is bounded and \ $\mu(K)<\infty$.

\noindent{(ii):} If \ $\mu$ \ is a measure on \ $K$ \ such that \ $\mu(K)<\infty$ \ and
 it is absolutely continuous with respect to the Lebesgue measure on \ $K$ \ with Radon-Nikodym derivative \ $h_\mu$, \ then
  \[
   F_{K,\mu}(x,y)=\int_K d_1((x,y),(\alpha,\beta))h_\mu(\alpha,\beta)\,\dd\alpha\dd\beta,
     \qquad (x,y)\in\RR^2.
 \]
With
 \[
   h_\mu(\alpha,\beta)
               :=\begin{cases}
                   \frac{1}{A(K)} & \text{if \ $(\alpha,\beta)\in K$,}\\
                   0 & \text{if \ $(\alpha,\beta)\not\in K$,}
                 \end{cases}
 \]
 we have \ $F_{K,\mu}$ \ coincides with \ $F_K$ \ given in Definition \ref{Def_gen_conic_VN}.
Note also that the conic function \ $F_K$ \ can be interpreted as the expectation of
 an appropriate random variable.
Namely, \ $F_K(x,y) = \EE [d_1((x,y),(\xi,\eta))]$, $(x,y)\in\RR^2$,
 \ where \ $(\xi,\eta)$ \ is a uniformly distributed random variable on \ $K$.
\proofend 
\end{Rem}

Next we generalize Theorems 3, 4 and 5, Lemmas 6 and 7 and Corollary 1 in Vincze and Nagy \cite{VinNag}
 for the generalized conic function \ $F_{K,\mu}$.

\begin{Thm}\label{Thm2}
The generalized conic function \ $F_{K,\mu}:\RR^2\to\RR_+$ \ is a convex function which satisfies the growth condition
 \begin{align*}
  \liminf_{\Vert (x,y)\Vert\to\infty} \frac{F_{K,\mu}(x,y)}{\sqrt{x^2+y^2}} \geq \mu(K)>0.
 \end{align*}
Consequently, the level sets of the function \ $F_{K,\mu}$ \ are bounded and hence compact subsets of \ $\RR^2$.
\end{Thm}

\noindent{\bf Proof.}
Recall that
\[
   F_{K,\mu}(x,y)=\int_K d_1((x,y),(\alpha,\beta))\,\mu(\dd\alpha,\dd\beta),
     \qquad (x,y)\in\RR^2.
 \]
The convexity of \ $F_{K,\mu}$ \ is clear, since the integrand is a convex function for any fixed element \ $(\alpha,\beta)\in K$,
 \ and the Lebesgue integral with respect to the measure \ $\mu$ \ is monotone.
Further, since \ $d_2((x,y),(\alpha,\beta))\leq d_1((x,y),(\alpha,\beta))$, $(x,y),(\alpha,\beta)\in\RR^2$, \
 where \ $d_2$ \ is the standard Euclidean distance on \ $\RR^2$, \ we have
 \[
    F_{K,\mu}(x,y) \geq \int_K d_2((x,y),(\alpha,\beta))\,\mu(\dd\alpha,\dd\beta),\qquad (x,y)\in\RR^2,
 \]
 and then
 \[
    \frac{F_{K,\mu}(x,y)}{\sqrt{x^2 + y^2}}
        \geq \int_K \left( \frac{d_2((x,y),(\alpha,\beta)) - \sqrt{x^2 + y^2}}{\sqrt{x^2 + y^2}} + 1\right)
                 \,\mu(\dd\alpha,\dd\beta)
 \]
 for \ $(x,y)\in\RR^2$, $(x,y)\ne(0,0)$.
The triangle inequality shows that
 \begin{align*}
   \sqrt{x^2 + y^2} = d_2((x,y),(0,0))
      &\leq d_2((x,y),(\alpha,\beta)) + d_2((\alpha,\beta),(0,0))\\
      & =  d_2((x,y),(\alpha,\beta)) + \sqrt{\alpha^2 + \beta^2},
 \end{align*}
 and then
 \begin{align*}
  \frac{F_{K,\mu}(x,y)}{\sqrt{x^2 + y^2}}
     \geq  \int_K \left( 1 - \frac{\sqrt{\alpha^2 + \beta^2}}{\sqrt{x^2 + y^2}} \right)
                 \,\mu(\dd\alpha,\dd\beta),
                 \qquad (x,y)\in\RR^2,\; (x,y)\ne(0,0).
 \end{align*}
By Fatou's lemma,
 \begin{align*}
    \liminf_{\Vert (x,y)\Vert\to\infty} \frac{F_{K,\mu}(x,y)}{\sqrt{x^2+y^2}}
      & \geq \liminf_{\Vert (x,y)\Vert\to\infty}
             \int_K \left( 1 - \frac{\sqrt{\alpha^2 + \beta^2}}{\sqrt{x^2 + y^2}} \right)
                 \,\mu(\dd\alpha,\dd\beta)\\
      & \geq \int_K \liminf_{\Vert (x,y)\Vert\to\infty}
                   \left( 1 - \frac{\sqrt{\alpha^2 + \beta^2}}{\sqrt{x^2 + y^2}} \right)
                 \,\mu(\dd\alpha,\dd\beta)
       = \mu(K)>0.
 \end{align*}
Here for completeness we note that one can use Fatou's lemma, since for all \ $c>0$,
 \begin{align*}
  &\int_K\inf\left\{ 1- \frac{\sqrt{\alpha^2+\beta^2}}{\sqrt{x^2+y^2}} : \Vert (x,y)\Vert\geq c\right\}
         \,\mu(\dd\alpha,\dd\beta)\\
  &\qquad\qquad = \int_K\left( 1 - \frac{\sqrt{\alpha^2+\beta^2}}{c}\right) \,\mu(\dd\alpha,\dd\beta)
                >-\infty,
 \end{align*}
 where the last inequality follows by that \ $K$ \ is compact (hence bounded) and \ $\mu(K)<\infty$.

Let \ $d\in\RR_+$ \ and let us suppose that the level set \ $\{(x,y)\in\RR^2 : F_{K,\mu}(x,y)\leq d\}$ \ is unbounded.
Then one can choose a sequence \ $(x_n,y_n)$, $n\in\NN$, \ such that \ $F_{K,\mu}(x_n,y_n)\leq d$, $n\in\NN$, \ and
 \ $\lim_{n\to\infty}\Vert (x_n,y_n)\Vert = \infty$.
\ This would imply that
 \[
     \lim_{n\to\infty} \frac{F_{K,\mu}(x_n,y_n)}{\sqrt{x_n^2+y_n^2}} = 0,
 \]
 which contradicts to the growth condition.
\proofend

\begin{Lem}\label{Lem6}
Let us suppose that Condition (C.3) holds.
For the generalized conic function \ $F_{K,\mu}$, \ we have
 \begin{align*}
  F_{K,\mu}(x,y)
     & = x\big( \mu(\{K<_1x\}) - \mu(\{x<_1 K\}) \big)
        - \int_K \alpha(\mathbf 1_{\{\alpha<x\}}  - \mathbf 1_{\{x< \alpha\}} )\,\mu(\dd\alpha,\dd\beta) \\
     &\phantom{=\;} + y\big( \mu(\{K<_2y\}) - \mu(\{y<_2 K\}) \big)
        - \int_K\!\beta(\mathbf 1_{\{\beta<y\}}  - \mathbf 1_{\{y<\beta\}} )\,\mu(\dd\alpha,\dd\beta)
 \end{align*}
 for all \ $(x,y)\in\RR^2$.
\end{Lem}

\noindent{\bf Proof.}
By definition,
 \[
   F_{K,\mu}(x,y) = \int_K(\vert x-\alpha\vert + \vert y-\beta\vert)\,\mu(\dd\alpha,\dd\beta),
      \qquad (x,y)\in\RR^2.
 \]
Here
 \begin{align*}
  &\int_K \vert x  - \alpha\vert \,\mu(\dd\alpha,\dd\beta)
     =   \int_{K<_1x} \vert x - \alpha\vert \,\mu(\dd\alpha,\dd\beta)
        +   \int_{x\leq_1 K} \vert x - \alpha\vert \,\mu(\dd\alpha,\dd\beta)\\
  & =   \int_{K<_1x} (x - \alpha) \,\mu(\dd\alpha,\dd\beta)
        +   \int_{x\leq_1 K} (\alpha - x) \,\mu(\dd\alpha,\dd\beta) \\
  & = x\big( \mu(\{K<_1x\}) - \mu(\{x\leq_1 K\}) \big)
      - \int_{K<_1x} \alpha \,\mu(\dd\alpha,\dd\beta)
      + \int_{x\leq_1 K} \alpha \,\mu(\dd\alpha,\dd\beta),
 \end{align*}
 and the integral \ $\int_K \vert y - \beta\vert \,\mu(\dd\alpha,\dd\beta)$ \
 can be handled similarly.
The assertion follows by taking into account Condition (C.3).
\proofend

\begin{Lem}\label{Lem5}
Let us suppose that Condition (C.3) holds.
For the generalized conic function \ $F_{K,\mu}$, \ we have
 \begin{align*}
  &D_1 F_{K,\mu}(x,y) =  \mu(\{K<_1x\}) - \mu(\{x<_1K\}), \qquad (x,y)\in\RR^2,\\
  &D_2 F_{K,\mu}(x,y) =  \mu(\{K<_2y\}) - \mu(\{y<_2K\}), \qquad (x,y)\in\RR^2.
 \end{align*}
\end{Lem}

\noindent{\bf Proof.}
Let \ $h>0$.
\ Then for all \ $(x,y)\in\RR^2$,
 \begin{align*}
  &\frac{F_{K,\mu}(x+h,y)  - F_{K,\mu}(x,y)}{h}
     = \int_K \frac{ \vert x+h -\alpha\vert - \vert x - \alpha \vert}{h} \,\mu(\dd\alpha,\dd\beta)\\
  &\quad\quad = \int_{K<_1x} \frac{ \vert x+h -\alpha\vert - \vert x - \alpha \vert}{h} \,\mu(\dd\alpha,\dd\beta)\\
  &\quad\quad\phantom{=\;}
       + \int_{x\leq_1 K\leq_1 x+h} \frac{ \vert x+h -\alpha\vert - \vert x - \alpha \vert}{h} \,\mu(\dd\alpha,\dd\beta)\\
  &\quad\quad\phantom{=\;}
       + \int_{x+h <_1 K} \frac{ \vert x+h -\alpha\vert - \vert x - \alpha \vert}{h} \,\mu(\dd\alpha,\dd\beta)\\
  &\quad\quad = \int_{K<_1x} \frac{ x+h -\alpha  - (x - \alpha) }{h} \,\mu(\dd\alpha,\dd\beta)\\
  &\quad\quad\phantom{=\;}
       + \int_{x\leq_1 K\leq_1 x+h} \frac{ x+h -\alpha - (\alpha - x) }{h} \,\mu(\dd\alpha,\dd\beta)\\
  &\quad\quad \phantom{=\;}
       + \int_{x+h <_1 K} \frac{ \alpha - x-h  - ( \alpha - x) }{h} \,\mu(\dd\alpha,\dd\beta) \\
  &\quad\quad = \mu(\{K<_1x\})  - \mu(\{x+h <_1 K\})\\
  &\quad\quad \phantom{=\;}
       + \int_{x\leq_1 K\leq_1 x+h} \frac{ \vert x+h -\alpha \vert - \vert x- \alpha \vert }{h} \,\mu(\dd\alpha,\dd\beta).
 \end{align*}
Using that \ $\vert \vert a\vert - \vert b\vert \vert \leq \vert a - b \vert$, $a,b\in\RR$, \ for the integrand we have
 \begin{align*}
   \left\vert  \frac{ \vert x+h -\alpha \vert - \vert x- \alpha \vert }{h} \right\vert
      \leq \frac{1}{h} \vert x+h-\alpha - (x-\alpha)\vert
      = \frac{\vert h\vert}{h}
      =1,\qquad x,\alpha\in\RR,\; h>0,
 \end{align*}
 and hence, by dominated convergence theorem,
 \begin{align*}
   &\left\vert \int_{x\leq_1 K\leq_1 x+h} \frac{ \vert x+h -\alpha \vert - \vert x- \alpha \vert }{h} \,\mu(\dd\alpha,\dd\beta) \right\vert\\
   &\qquad\leq \int_{x\leq_1 K\leq_1 x+h} \left\vert \frac{ \vert x+h -\alpha \vert - \vert x- \alpha \vert }{h} \right\vert
       \,\mu(\dd\alpha,\dd\beta) \\
   &\qquad\leq \mu(\{ x\leq_1 K\leq_1 x+h\})
    \to \mu(\{K=_1x\})=0
 \end{align*}
 as \ $h\downarrow 0$.
\ Then, for all \ $(x,y)\in\RR^2$,
 \begin{align}\label{seged12}
  \begin{split}
   \lim_{h\downarrow 0} \frac{F_{K,\mu}(x+h,y)  - F_{K,\mu}(x,y)}{h}
   &  = \mu(\{K<_1x\}) - \mu(\{x\leq_1 K\})\\
   & = \mu(\{K<_1x\}) - \mu(\{x<_1K\}).
  \end{split}
 \end{align}

Similarly, if \ $h<0$, \ then
 \begin{align*}
    \frac{F_{K,\mu}(x+h,y)  - F_{K,\mu}(x,y)}{h}
      &= \mu(\{K<_1x+h\}) - \mu(\{x <_1 K\}) \\
      &\phantom{=\;}  + \int_{x+h\leq_1 K\leq_1 x}
                        \frac{\vert x + h - \alpha \vert - \vert x -\alpha \vert }{h} \,\mu(\dd\alpha,\dd\beta)
 \end{align*}
 for all \ $(x,y)\in\RR^2$, \ and hence, using again Condition (C.3),
 \begin{align}\label{seged13}
  \begin{split}
  \lim_{h\uparrow 0} \frac{F_{K,\mu}(x+h,y)  - F_{K,\mu}(x,y)}{h}
   & = \mu(\{K\leq_1 x\}) - \mu(\{x<_1 K\})\\
   &\phantom{=\;} = \mu(\{K<_1x\}) - \mu(\{x<_1K\})
  \end{split}
 \end{align}
 for all \ $(x,y)\in\RR^2$.
\ Then \eqref{seged12} and \eqref{seged13} yield that \ $D_1 F_{K,\mu}(x,y) = \mu(\{K<_1x\}) - \mu(\{x<_1K\})$, $(x,y)\in\RR^2$.

In a similar way, we have \ $D_2 F_{K,\mu}(x,y) = \mu(\{K<_2y\}) - \mu(\{y<_2K\})$, $(x,y)\in\RR^2$.
\proofend

If \ $\mu$ \ is a measure on \ $(\RR^d,\cB(\RR^d))$, \ then by the \ $\mu$-area of a Borel measurable set \ $S\in\cB(\RR^d)$,
 \ we mean \ $\mu(S)$.

\begin{Cor}\label{Cor1}
Let us suppose that Condition (C.3) holds.
A point in \ $\RR^2$ \ is a global minimizer of the generalized conic function \ $F_{K,\mu}$ \ if and only if
 it bisects the \ $\mu$-area of \ $K$, \ i.e., the vertical and the horizontal lines through this point cut the body
 \ $K$ \ into two parts with equal \ $\mu$-areas.
Moreover, if Conditions (C.1) and (C.2) hold too, then the convex function \ $F_{K,\mu}$ \ has a unique global minimizer
 \ $(x^*,y^*)\in\RR^2$, \ that is, \ $F_{K,\mu}(x,y)>F_{K,\mu}(x^*,y^*)$ \ for \ $(x,y)\ne(x^*,y^*)$, \ $(x,y)\in\RR^2$.
\end{Cor}

\noindent{\bf Proof.}
First note that under Condition (C.3) the concept of bisection of the \ $\mu$-area of \ $K$ \ is well-defined.
The first part of the corollary is a consequence of Lemma \ref{Lem5} using that a local minimum of a convex function defined
 on \ $\RR^2$ \ is a global minimum, too.
Under Conditions (C.1), (C.2) and (C.3), the existence of a global minimizer \ $(x^*,y^*)$ \ of \ $F_{K,\mu}$ \
 follows by that \ $F_{K,\mu}$ \ is a convex function defined on \ $\RR^2$ \ and its level sets are compact subsets of
 \ $\RR^2$ \ (see Theorem \ref{Thm2}).
Indeed, a finite-valued convex function defined on \ $\RR^2$ \ is continuous and it reaches its minimum on every compact set.
Now we turn to prove the uniqueness of \ $(x^*,y^*)$.
\ The proof goes along the very same lines as in the proof of Proposition \ref{Pro1}.
Indeed, the area \ $A$ \ ($2$-dimensional Lebesgue measure) has to be replaced by the measure \ $\mu$.
\proofend

Before we generalize Theorem 4 in Vincze and Nagy \cite{VinNag} we need to introduce some notations and to recall the Cavalieri principle
 for product measures.

\begin{Def}\label{Def_product_measure}
Let \ $\mu_1$ \ and \ $\mu_2$ \ be \ $\sigma$-finite measures on \ $(\RR,\cB(\RR))$ \ and let
 \ $\mu:=\mu_1\times\mu_2$ \ be their product measure on \ $(\RR^2,\cB(\RR^2))$.
\ Given a measurable set \ $S\in\cB(\RR^2)$, \ the generalized \ $X$-ray functions of \ $S$ \ with respect to \ $\mu$ \
 into the coordinate directions are defined by
 \[
    X_{S,\mu}(y):= \mu_1(S_y),\quad y\in\RR,\qquad \text{and}\qquad   Y_{S,\mu}(x):=\mu_2(S_x),\quad x\in\RR,
 \]
 where \ $S_x:=\{ y\in\RR : (x,y)\in S\}$ \ and \ $S_y:=\{ x\in\RR : (x,y)\in S\}$.
\ (Note that \ $S_x,S_y\in\cB(\RR)$ \ for all \ $x,y\in\RR$, \ see, e.g., Lemma 5.1.1 in Cohn \cite{Coh}.)
\end{Def}

For the product measure \ $\mu$ \ defined in Definition \ref{Def_product_measure}, we have \ $\mu(K)<\infty$.

\begin{Thm}\label{Thm_Cavalieri}{\bf(The Cavalieri principle, see, e.g., Cohn \cite[Theorem 5.1.3]{Coh})}
Let \ $\mu_1$ \ and \ $\mu_2$ \ be \ $\sigma$-finite measures on \ $(\RR,\cB(\RR))$ \ and let
 \ $\mu:=\mu_1\times\mu_2$ \ be their product measure on \ $(\RR^2,\cB(\RR^2))$.
\ If \ $S\in\cB(\RR^2)$, \ then the functions \ $X_{S,\mu}, Y_{S,\mu} : \RR\to\RR_+$ \ are Borel measurable, and
 \[
   \mu(S) = (\mu_1\times\mu_2)(S) = \int_\RR Y_{S,\mu}(x)\mu_1(\dd x) =  \int_\RR X_{S,\mu}(y)\mu_2(\dd y).
 \]
\end{Thm}

\begin{Thm}\label{Thm1}
Let \ $K,K^*\subset \RR^2$ \ be compact bodies,
 let \ $\mu_i$, $\mu_i^*$, $i=1,2$, \ be \ $\sigma$-finite measures on \ $(\RR,\cB(\RR))$ \
 that are absolutely continuous with respect to the Lebesgue measure on \ $(\RR,\cB(\RR))$ \
 with Radon-Nikodym derivatives \ $f_i$, $f_i^*$, $i=1,2$.
\ Let \ $\mu:=\mu_1\times\mu_2$ \ and \ $\mu^*:=\mu_1^*\times\mu_2^*$ \ be their product measures
 on \ $(\RR^2,\cB(\RR^2))$ \  and we assume that \ $\mu$ \ and \ $\mu^*$ \ are supported by \ $K$ \ and
 \ $K^*$, \ respectively.
Let us suppose that Condition (C.3) holds for \ $K$ \ and \ $\mu$, \ and \ $K^*$ \ and \ $\mu^*$, \ respectively.
Then \ $F_{K,\mu} = F_{K^*,\mu^*}$ \ if and only if \ $f_2(y)X_{K,\mu}(y) = f_2^*(y)X_{K^*,\mu^*}(y)$ \ for
 (Lebesgue) almost every \ $y\in\RR$, \ and \ $f_1(x)Y_{K,\mu}(x) = f_1^*(x)Y_{K^*,\mu^*}(x)$ \
 for (Lebesgue) almost every \ $x\in\RR$.
\end{Thm}

\noindent{\bf Proof.}
By Theorem \ref{Thm_Cavalieri} (the Cavalieri principle), for all \ $x,y\in\RR$,
 \begin{align}\label{seged14}
   \begin{split}
   &\mu(K<_1x) = \int_\RR Y_{K<_1x,\mu}(s)\,\mu_1(\dd s)
               = \int_{-\infty}^x Y_{K,\mu}(s)\,\mu_1(\dd s)
               = \int_{-\infty}^x Y_{K,\mu}(s)f_1(s)\,\dd s,\\
   &\mu(x<_1K) = \int_\RR Y_{x<_1K,\mu}(s)\,\mu_1(\dd s)
               = \int_x^{\infty} Y_{K,\mu}(s)\,\mu_1(\dd s)
               = \int_x^{\infty} Y_{K,\mu}(s)f_1(s)\,\dd s,\\
   &\mu(K<_2y) = \int_\RR X_{K<_2y,\mu}(t)\,\mu_2(\dd t)
               = \int_{-\infty}^y X_{K,\mu}(t)\,\mu_2(\dd t)
               = \int_{-\infty}^y X_{K,\mu}(t)f_2(t)\,\dd t,\\
   & \mu(y<_2K) = \int_\RR X_{y<_2K,\mu}(t)\,\mu_2(\dd t)
                = \int_y^{\infty} X_{K,\mu}(t)\,\mu_2(\dd t)
                = \int_y^{\infty} X_{K,\mu}(t)f_2(t)\,\dd t,
  \end{split}
 \end{align}
 and, by Fubini's theorem, for all \ $x,y\in\RR$,
 \begin{align}\label{seged15}
   \begin{split}
   & \int_K \alpha\mathbf 1_{\{\alpha<x\}}\,\mu(\dd\alpha,\dd\beta)
       = \int_{-\infty}^x s Y_{K,\mu}(s)\,\mu_1(\dd s)
       = \int_{-\infty}^x s Y_{K,\mu}(s)f_1(s)\,\dd s,\\
   & \int_K \alpha\mathbf 1_{\{x<\alpha\}}\,\mu(\dd\alpha,\dd\beta)
       = \int_x^{\infty} s Y_{K,\mu}(s)\,\mu_1(\dd s)
       = \int_x^{\infty} s Y_{K,\mu}(s)f_1(s)\,\dd s,\\
   & \int_K \beta\mathbf 1_{\{\beta<y\}}\,\mu(\dd\alpha,\dd\beta)
       = \int_{-\infty}^y t X_{K,\mu}(t)\,\mu_2(\dd t)
       = \int_{-\infty}^y t X_{K,\mu}(t)f_2(t)\,\dd t,\\
   & \int_K \beta\mathbf 1_{\{y<\beta\}}\,\mu(\dd\alpha,\dd\beta)
       = \int_y^{\infty} t X_{K,\mu}(t)\,\mu_2(\dd t)
       = \int_y^{\infty} t X_{K,\mu}(t)f_2(t)\,\dd t.
  \end{split}
 \end{align}
Indeed, for example, the first statement of \eqref{seged15} holds since, by Fubini's theorem for
 non-rectangular regions,
 \begin{align*}
  \int_K & \alpha\mathbf 1_{\{\alpha<x\}}\,\mu(\dd\alpha,\dd\beta)
      = \int_{\alpha_b}^{\alpha_u} \left( \int_{K_\alpha}\alpha\mathbf 1_{\{\alpha<x\}}\,\mu_2(\dd\beta)\right)\mu_1(\dd\alpha)\\
    & = \int_{\alpha_b}^{\alpha_u} \alpha\mathbf 1_{\{\alpha<x\}} \mu_2(K_\alpha) \,\mu_1(\dd\alpha)
      = \int_{\alpha_b}^{\alpha_u} \alpha\mathbf 1_{\{\alpha<x\}} Y_{K,\mu}(\alpha) \,\mu_1(\dd\alpha)\\
     & = \int_{-\infty}^x s Y_{K,\mu}(s) \,\mu_1(\dd s),
 \end{align*}
 where \ $K_\alpha = \{ \beta \in\RR \mid (\alpha,\beta)\in K\}$ \ and
 \begin{align*}
  \alpha_b:=\inf\big\{ \alpha \mid \exists\,\beta\in\RR : (\alpha,\beta)\in K\big\},\quad
  \alpha_u:=\sup\big\{ \alpha \mid \exists\,\beta\in\RR : (\alpha,\beta)\in K\big\}.
 \end{align*}
Further, by \eqref{seged14}, Lemma \ref{Lem5} and Lebesgue differentiation theorem,
 \begin{align}\label{seged16}
  \begin{split}
   D_1D_1 F_{K,\mu}(x,y)
     & = D_1\big(\mu(\{K<_1x\}) - \mu(\{x<_1 K\})\big) \\
     & = D_1\left( \int_{-\infty}^x Y_{K,\mu}(s)f_1(s)\,\dd s
                   - \int_x^{\infty} Y_{K,\mu}(s)f_1(s)\,\dd s \right) \\
     & = 2Y_{K,\mu}(x)f_1(x) \quad \text{for all \ $y\in\RR$ \ and almost every \ $x\in\RR$,}
  \end{split}
 \end{align}
 and, similarly,
 \begin{align}\label{seged16_kieg}
  \begin{split}
  &D_1D_2 F_{K,\mu}(x,y) = D_2D_1 F_{K,\mu}(x,y) = 0 \quad \text{for all \ $(x,y)\in\RR^2$,}\\
  &D_2D_2 F_{K,\mu}(x,y) = 2X_{K,\mu}(y)f_2(y) \quad \text{for all \ $x\in\RR$ \ and almost every \ $y\in\RR$.}
  \end{split}
 \end{align}

Let us suppose that \ $F_{K,\mu} = F_{K^*,\mu^*}$.
\ By \eqref{seged16} and \eqref{seged16_kieg}, we have
 \ $f_1(x)Y_{K,\mu}(x) = f_1^*(x)Y_{K^*,\mu^*}(x)$ \ for almost every \ $x\in\RR$, \
 and \ $f_2(y)X_{K,\mu}(y) = f_2^*(y)X_{K^*,\mu^*}(y)$ \ for almost every \ $y\in\RR$, \ as desired.

Conversely, let us suppose that \ $f_2(y)X_{K,\mu}(y) = f_2^*(y)X_{K^*,\mu^*}(y)$ \ for almost every \ $y\in\RR$, \ and
 \ $f_1(x)Y_{K,\mu}(x) = f_1^*(x)Y_{K^*,\mu^*}(x)$ \ for almost every \ $x\in\RR$.
\ Then, by Lemma \ref{Lem6}, \eqref{seged14} and \eqref{seged15}, we get \ $F_{K,\mu} = F_{K^*,\mu^*}$.
\proofend

\begin{Rem}
Note that, under the conditions of Theorem \ref{Thm1}, for almost every \ $(x,y)\in\RR^2$, \
 the matrix consisting of the second order partial derivatives of \ $F_{K,\mu}$ \ takes the form
 \[
   \begin{bmatrix}
     2f_1(x)Y_{K,\mu}(x) & 0 \\
     0 & 2f_2(y)X_{K,\mu}(y) \\
   \end{bmatrix},
 \]
 which is a positive semidefinite matrix, since the Radon-Nikodym derivatives \ $f_i$ \ and \ $f_i^*$, $i=1,2$,
 \ are non-negative almost everywhere.
Note also that this is in accordance with the fact that \ $F_{K,\mu}$ \ is a convex function due to Theorem \ref{Thm2}.
\proofend
\end{Rem}

Before we generalize Theorem 5 in Vincze and Nagy \cite{VinNag}, we need to recall some notions.

\begin{Def}
Let \ $K$ \ be a compact body in \ $\RR^2$.
\ For all \ $\varepsilon>0$, \ the outer parallel body \ $K^\varepsilon$ \ is the union of closed Euclidean balls centered at the points
 of \ $K$ \ with radius \ $\varepsilon>0$.
\end{Def}

\begin{Def}
The Hausdorff distance between two compact bodies \ $K$ \ and \ $L$ \ is given by
 \[
   H(K,L):=\inf\big\{\varepsilon>0 : K\subset L^\varepsilon \quad \text{and}\quad L\subset K^\varepsilon\big\}.
 \]
\end{Def}

The collection of compact bodies in \ $\RR^2$ \ furnished with the Hausdorff distance \ $H$ \ is a metric space,
 see, e.g., Beer \cite{Bee}.

\begin{Lem}\label{Lemma_regular}
Let \ $K_n$, $n\in\NN$, \ $K$ \ be compact bodies, and let \ $\mu$ \ be a Radon measure on \ $(\RR^2,\cB(\RR^2))$.
\renewcommand{\labelenumi}{{\rm(\roman{enumi})}}
 \begin{enumerate}
  \item We have \ $\lim_{\varepsilon\downarrow 0} \mu(K^\varepsilon) = \mu(K)$.
  \item If \ $K_n\to K$ \ as \ $n\to\infty$ \ with respect to the Hausdorff metric \ $H$,
             \ then the following regularity properties are equivalent:
             \begin{enumerate}
               \item[(a)] $\lim_{n\to\infty} \mu((K\setminus K_n) \cup (K_n\setminus K)) = 0$,
               \item[(b)] $\lim_{n\to\infty} \mu(K_n) = \mu(K)$.
             \end{enumerate}
\end{enumerate}
\end{Lem}

\noindent{\bf Proof.}
The proofs go along the very same lines as those of Lemmas 1 and 2 in Vincze and Nagy \cite{VinNag} by replacing
 the area \ $A$ \ ($2$-dimensional Lebesgue measure) by the measure \ $\mu$ \ in the proofs and
 refereeing to that \ $\mu(L)<\infty$ \ for all compact sets \ $L\subset \RR^2$ \
 (due to that \ $\mu$ \ is a Radon measure).
\proofend

\begin{Def}
Let \ $K_n$, $n\in\NN$, \ and \ $K$ \ be compact bodies, and let \ $\mu$ \ be a Radon measure on \ $(\RR^2,\cB(\RR^2))$.
The convergence \ $K_n\to K$ \ as \ $n\to\infty$ \ with respect to the Hausdorff metric is called regular if one of
 the conditions (a) and (b) of part (ii) of Lemma \ref{Lemma_regular} holds.
\end{Def}

\begin{Thm}
Let \ $K_n$, $n\in\NN$, \ and \ $K$ \ be compact bodies, and let \ $\mu$ \ be a Radon measure on \ $(\RR^2,\cB(\RR^2))$
 \ supported by \ $K^\varepsilon$ \ for some \ $\varepsilon>0$.
\ Let us suppose that the convergence \ $K_n\to K$ \ as \ $n\to\infty$ \ with respect to the Hausdorff metric is regular.
Then
 \[
  \lim_{n\to\infty} F_{K_n,\mu}(x,y) = F_{K,\mu}(x,y),\qquad (x,y)\in\RR^2.
 \]
\end{Thm}

\noindent{\bf Proof.}
The proof goes along the very same lines as that of Theorem 5 in Vincze and Nagy \cite{VinNag}, but replacing the integration
 with respect to the two-dimensional Lebesgue measure by the integration with respect to the measure \ $\mu$.
\proofend

For the remaining sections of the paper we will need some further properties of the convex function \ $F_{K,\mu}$.
\ Next we recall some general facts from the theory of convex functions, see, e.g., Polyak \cite[Lemma 3, Section 1.1.4]{Pol}.

\begin{Lem}\label{Lem7}
Let \ $F:\RR^d\to\RR$ \ be a differentiable and convex function such that its gradient is Lipschitz continuous with constant
 \ $L>0$, \ i.e.
 \begin{align}\label{seged19}
     \|\mathrm{grad}\ F(p)-\mathrm{grad}\ F (q)\|\leq L \|p-q\|,\qquad p,q\in\RR^d,
 \end{align}
 where \ $\mathrm{grad}\ F(p):= (D_1F(p),D_2F(p))^\top$, $p\in\RR^d$.
\ Then we have an affine lower bound
 \[
 F(q)\geq F(p)+\langle \mathrm{grad}\ F(p), q-p \rangle,\qquad p,q\in\RR^d.
 \]
\end{Lem}

\begin{Lem}\label{Lem_conic_uniform}
Let \ $\mu_1$ \ and \ $\mu_2$ \ be \ $\sigma$-finite measures on \ $(\RR,\cB(\RR))$ \ that are absolutely continuous
 with respect to the Lebesgue measure on \ $(\RR,\cB(\RR))$ \ with bounded Radon-Nikodym derivatives.
Let \ $\mu:=\mu_1\times\mu_2$ \ be their product measure on \ $(\RR^2,\cB(\RR^2))$ \ and we assume that \ $\mu$ \ is
 supported by \ $K$.
\ Further, let us suppose that Condition (C.3) holds.
Then the generalized conic function \ $F_{K,\mu}:\RR^2\to\RR$ \ associated with K and \ $\mu$ \ satisfies the
 conditions of Lemma \ref{Lem7}, and, consequently, we have an affine lower bound for \ $F_{K,\mu}$.
\end{Lem}

\noindent{\bf Proof.}
By Theorem \ref{Thm2}, \ $F_{K,\mu}$ \ is convex.
Under Condition (C.3), by Lemma \ref{Lem5} and \eqref{seged14},
 \begin{align*}
  D_1 F_{K,\mu}(x,y)
    & = \int_{-\infty}^x Y_{K,\mu}(s)\,\mu_1(\dd s)
       - \int_x^{\infty} Y_{K,\mu}(s)\,\mu_1(\dd s) \\
    & = \int_{-\infty}^x Y_{K,\mu}(s)f_1(s)\,\mu_1(\dd s)
       - \int_x^{\infty} Y_{K,\mu}(s)f_1(s)\,\mu_1(\dd s)
 \end{align*}
 for \ $(x,y)\in\RR^2$, \ where \ $f_1$ \ denotes the (bounded)
 Radon-Nikodym derivative of \ $\mu_1$ \ with respect to the Lebesgue measure on \ $\RR$.
\ Using that the integral as a function of the upper limit of the integration is continuous, we have \ $D_1 F_{K,\mu}$ \
 is continuous on \ $\RR^2$.
\ Similarly, one can check that \ $D_2 F_{K,\mu}$ \ is also continuous on \ $\RR^2$.
\ This implies that \ $F_{K,\mu}$ \ is differentiable on \ $\RR^2$.

Condition \eqref{seged19} for \ $F_{K,\mu}$ \ can be checked as follows.
Let us start with the difference of the partial derivatives with respect to the first variable
 \begin{align*}
 D_1F_{K,\mu}(q)&-D_1F_{K,\mu}(p)\\
                & =\mu(K<_1q^{(1)})-\mu(q^{(1)}<_1K)-(\mu(K<_1p^{(1)})-\mu(p^{(1)}<_1K))
 \end{align*}
 for all \ $p=(p^{(1)},p^{(2)}),\, q=(q^{(1)},q^{(2)})\in \RR^2$, \ where the equality follows by Lemma \ref{Lem5}.
\ We have
 $$
 \mu(K<_1q^{(1)})=\mu(K<_1\min \{p^{(1)},q^{(1)}\})+\mu(\min\{p^{(1)},q^{(1)}\}<_1K<_1q^{(1)})
 $$
 and
 $$
 \mu(q^{(1)}<_1K)=\mu(\max \{p^{(1)},q^{(1)}\}<_1K)+\mu(q^{(1)}<_1K<_1\max \{p^{(1)},q^{(1)}\}).
 $$
Of course we can change the role of \ $q$ \ and \ $p$ \ to express
 \ $\mu(K<_1p^{(1)})$ \ and \ $\mu(p^{(1)}<_1K)$ \ in a similar way.
Then
\begin{align*}
 D_1F_{K,\mu}(q)&-D_1F_{K,\mu}(p)\\
  &= \mu(\min\{p^{(1)},q^{(1)}\}<_1K<_1q^{(1)})-\mu(q^{(1)}<_1K<_1\max \{p^{(1)},q^{(1)}\}) \\
  &\phantom{=\;}- \mu(\min\{p^{(1)},q^{(1)}\}<_1K<_1p^{(1)})+\mu(p^{(1)}<_1K<_1\max \{p^{(1)},q^{(1)}\}).
 \end{align*}
Hence we can see that if $p^{(1)}=\min\{p^{(1)},q^{(1)}\}$ and,
consequently, $q^{(1)}=\max\{p^{(1)},q^{(1)}\}$, then
\[
 D_1F_{K,\mu}(q)-D_1F_{K,\mu}(p)=2\mu(p^{(1)}<_1K<_1 q^{(1)}).
 \]
If \ $q^{(1)}=\min\{p^{(1)},q^{(1)}\}$ \ and \ $p^{(1)}=\max\{p^{(1)},q^{(1)}\}$, \ then
 \[
  D_1F_{K,\mu}(q)-D_1F_{K,\mu}(p)=-2\mu(q^{(1)}<_1K<_1p^{(1)}).
  \]
In general,
 \[
  |D_1F_{K,\mu}(q)-D_1F_{K,\mu}(p)|=2\mu(\min \{p^{(1)},q^{(1)}\}<_1 K<_1 \max
   \{p^{(1)},q^{(1)}\}).
 \]
Therefore, using Theorem \ref{Thm_Cavalieri} (the Cavalieri principle), we can estimate the difference of the
absolute value of the first order partial derivatives of \ $F_{K,\mu}$ \ as follows
 \begin{align*}
   |D_1F_{K,\mu}(q)-D_1F_{K,\mu}(p)|
   &\leq 2 \int_{\min \{p^{(1)},q^{(1)}\}}^{\max \{p^{(1)},q^{(1)}\}} Y_{K,\mu}(s)\,\mu_1(\dd s)\\
   &\leq 2 \left(\sup_{s\in\RR} Y_{K,\mu}(s)\right)
           \mu_1\big( \big(\min \{p^{(1)},q^{(1)}\},\max \{p^{(1)},q^{(1)}\}\big)\big)\\
   &= 2 \left(\sup_{s\in\RR} Y_{K,\mu}(s)\right)
       \int_{\min \{p^{(1)},q^{(1)}\}}^{\max \{p^{(1)},q^{(1)}\}} f_1(s)\,\dd s\\
   & \leq 2C_1 \left(\sup_{s\in\RR} Y_{K,\mu}(s)\right)
      |p^{(1)}-q^{(1)}|
 \end{align*}
 with some constant \ $C_1>0$, \ where \ $\sup_{s\in\RR} Y_{K,\mu}(s)<\infty$ \ (since \ $\mu(K)<\infty$),
 \ and \ $f_1$ \ denotes the bounded Radon-Nikodym derivative of \ $\mu_1$ \ with respect to the Lebesgue measure on
  \ $\RR$.
\ Similarly,
 \[
  |D_2F_{K,\mu}(q)-D_2F_{K,\mu}(p)|\leq 2 C_2\left(\sup_{t\in\RR} X_{K,\mu}(t)\right) |p^{(2)}-q^{(2)}|
 \]
 with some constant \ $C_2>0$.
\ Therefore
 \begin{align*}
   &\Vert \mathrm{grad}\ F_{K,\mu}(p)-\mathrm{grad}\ F_{K,\mu}(q) \Vert\\
   &\qquad = \sqrt{(D_1F_{K,\mu}(p) - D_1F_{K,\mu}(q))^2 + (D_2F_{K,\mu}(p) - D_2F_{K,\mu}(q))^2}\\
   &\qquad\leq L \Vert p - q \Vert, \qquad p,q\in\RR^2,
 \end{align*}
 where
 \[
  L:=2\max\left\{C_1\sup_{s\in\RR} Y_{K,\mu}(s) ,C_2\sup_{t\in\RR} X_{K,\mu}(t) \right\},
 \]
 i.e., condition \eqref{seged19} for \ $F_{K,\mu}$ \  is satisfied with \ $d=2$ \ and with the Lipschitz constant \ $L$ \ given above.
\proofend

\section{A stochastic algorithm for the global minimizer of \ $F_{K,\mu}$}
\label{Section_stoch_alg_gen_conic}

We provide a stochastic algorithm for computing the global minimizer of generalized conic function \ $F_{K,\mu}$ \
 introduced in Definition \ref{Def_gen_conic}, and we prove almost sure and \ $L^q$-convergence of this algorithm.

In this section we assume that
 \[
   \mathbf{(C.4)} \qquad\qquad \text{$\mu$ \ is a probability measure on \ $K$.}
 \]

Let \ $(t_k)_{k\in\NN}$ \ be a decreasing sequence of positive numbers such that \ $\sum_{k=1}^\infty t_k=\infty$ \ and
 \ $\sum_{k=1}^\infty t_k^2<\infty$.

Let \ $(P_k)_{k\in\NN}$ \ be a sequence of independent identically distributed (2-dimensional) random variables
 such that their common distribution on \ $(\RR^2,\cB(\RR^2))$ \ is given by \ $\mu$.
\ Let \ $x_0\in K$ \ be arbitrarily chosen.
We define recursively a Markov chain \ $(X_k)_{k\in\ZZ_+}$ \ by
 \begin{align}\label{recursion}
   X_0:=x_0,\qquad \text{and} \qquad X_{k+1}:= X_k - t_{k+1}Q_{k+1}, \quad k\in\ZZ_+,
 \end{align}
 where
 \begin{align*}
   Q_{k+1} :=
      \begin{cases}
          \begin{pmatrix}
            1 \\
            1 \\
          \end{pmatrix}
            & \text{if \ $X_k^{(1)}\geq P_{k+1}^{(1)}$ \ and \ $X_k^{(2)}\geq P_{k+1}^{(2)}$,}\\
          \begin{pmatrix}
            1 \\
            -1 \\
          \end{pmatrix}
            & \text{if \ $X_k^{(1)}\geq P_{k+1}^{(1)}$ \ and \ $X_k^{(2)}< P_{k+1}^{(2)}$,}\\
          \begin{pmatrix}
            -1 \\
            1 \\
          \end{pmatrix}
            & \text{if \ $X_k^{(1)}<P_{k+1}^{(1)}$ \ and \ $X_k^{(2)}\geq P_{k+1}^{(2)}$,}\\
         \begin{pmatrix}
            -1 \\
            -1 \\
          \end{pmatrix}
            & \text{if \ $X_k^{(1)}< P_{k+1}^{(1)}$ \ and \ $X_k^{(2)}<P_{k+1}^{(2)}$, }
      \end{cases}
 \end{align*}
 with the notations \ $X_k:=(X_k^{(1)},X_k^{(2)})$, $P_k:=(P_k^{(1)},P_k^{(2)})$, $k\in\NN$.

\begin{Rem}
Note that if \ $\mu$ \ is a probability measure on \ $K$ \ such that it is absolutely continuous
 with respect to the Lebesgue measure on \ $K$ \ with Radon-Nikodym derivative (density function) \ $h_\mu$ \ given by
 \[
   h_\mu(x,y)=\begin{cases}
                   \frac{1}{A(K)} & \text{if \ $(x,y)\in K$,}\\
                   0 & \text{if \ $(x,y)\not\in K$,}
                 \end{cases}
 \]
 i.e., \ $\mu$ \ is the uniform distribution on \ $K$, \
 then \ $(P_k)_{k\in\NN}$ \ is a sequence of independent identically distributed (2-dimensional) random variables
 such that their common distribution is the uniform distribution on \ $K$.
 \proofend
\end{Rem}

\subsection{Almost sure and \ $L^q$-convergence of \ $(X_k)_{k\in\ZZ_+}$}

First we recall the so-called Robbins-Monro algorithm based on Bouleau and L\'epingle \cite[Theorem B.5.1, Chapter 2]{BouLep}.
This algorithm (in dimension 1) was originally invented by Robbins and Monro \cite{RobMon}.

Let \ $d\in\NN$ \ and \ $(t_n)_{n\in\ZZ_+}$ \ be a decreasing sequence of positive real numbers.
Let us suppose that all the random variables introduced below are defined on a probability space
 \ $(\Omega,\cF,\PP)$.
\ The Robbins-Monro algorithm generates a sequence of \ $\RR^d$-valued random variables
 \ $(\theta_n)_{n\in\ZZ_+}$ \ given by the recursion
 \begin{align*}
   \theta_{n+1} := \theta_n + t_{n+1}(\beta - \xi_{n+1}), \quad n\in\ZZ_+,
 \end{align*}
 where \ $\beta\in\RR^d$, \ $\theta_0$ \ is a given \ $\RR^d$-valued random variable,
 and \ $(\xi_n)_{n\in\ZZ_+}$ \ is a sequence of \ $d$-dimensional random variables
 such that there exists a Borel measurable function \ $M:\RR^d\to\RR^d$ \ satisfying
 \[
  \EE(\xi_{n+1} \mid \cF_n) = M(\theta_n) \qquad \text{$\PP$-almost surely for all \ $n\in\NN$,}
 \]
 where the filtration \ $(\cF_n)_{n\in\ZZ_+}$ \ is defined by \ $\cF_0:=\sigma(\theta_0)$ \
 (the sigma-algebra generated by \ $\theta_0$) \ and
 \ $\cF_n:=\sigma(\theta_0,\theta_1,\ldots,\theta_n,\xi_1,\ldots,\xi_n)$, \ $n\in\NN$
 \ (the sigma-algebra generated by \ $\theta_0,\theta_1,\ldots,\theta_n,\xi_1,\ldots,\xi_n$).

The following assumptions will be used.

\noindent {\bf Assumption (A.1)}: The \ $\RR^d$-valued random variable \ $\theta_0$ \ belongs to \ $L^q(\Omega,\cF,\PP)$,
 \ where \ $q\in\NN$.

\noindent {\bf Assumption (A.2)}: There exists some \ $B>0$ \ such that \ $\Vert\xi_n\Vert\leq B$ \ for all \ $n\in\NN$.

\noindent {\bf Assumption (A.3)}: There exists some \ $\theta^*\in\RR^d$ \ such that for each \ $\varepsilon\in(0,1)$,
 \[
     \inf_{\varepsilon\leq \Vert \theta - \theta^*\Vert\leq 1/\varepsilon}
            \langle \theta - \theta^*, M(\theta) -\beta \rangle >0,
 \]
 where \ $\langle \cdot,\cdot\rangle$ \ denotes the usual inner product in \ $\RR^d$.
Here Assumption (A.3) could be interpreted as a ''half-space'' assumption: roughly speaking,
 given the value of \ $\theta_n$, \ the expected value of \ $\theta_{n+1}$ \ will be
 on that side of the hyperplane through \ $\theta_n$ \ having normal vector \ $\theta^*-\theta_n$ \
 which contains \ $\theta^*$.

\begin{Thm}\label{Thm_Robbins_Monro}{\bf [Almost sure and \ $L^q$-convergence of Robbins-Monro algorithm]}
Let us suppose that Assumptions (A.1), (A.2) and (A.3) hold and that the decreasing sequence
 \ $(t_n)_{n\in\ZZ_+}$ \ of positive numbers satisfies
 \[
    \sum_{n=0}^\infty t_n = \infty \quad \text{and} \quad \sum_{n=0}^\infty t_n^2<\infty.
 \]
Then \ $\PP(\lim_{n\to\infty}\theta_n = \theta^*)=1$ \ and \ $\lim_{n\to\infty} \EE\Vert \theta_n - \theta^*\Vert^q=0$ \ for all
 \ $q\in\NN$.
\end{Thm}

Note that under the conditions of Theorem \ref{Thm_Robbins_Monro}
 the point \ $\theta^*\in\RR^d$ \ exists uniquely due to that, by Theorem \ref{Thm_Robbins_Monro},
 \ $\PP(\lim_{n\to\infty}\theta_n = \theta^*)=1$ \ and the limit of an almost surely convergent sequence of random variables
 is unique (up to probability one).
We also mention that, from a technical point of view, Assumption (A.3) is used for defining
 an appropriate non-negative supermartingale in order to prove the almost sure convergence
 of the sequence \ $(\theta_n)_{n\in\ZZ_+}$, \
 see, e.g., Bouleau and L\'epingle \cite[proof of Theorem B.5.1, Chapter 2]{BouLep}.

We will prove almost sure and \ $L^q$-convergence of the recursion given in \eqref{recursion}.
But first we present an auxiliary lemma.

\begin{Lem}
Let us consider the sequence \ $(X_k)_{k\in\ZZ_+}$ \ defined by \eqref{recursion}.
Let us suppose that Conditions (C.3) and (C.4) hold.
Then
 \begin{align}\label{seged8}
   \EE(Q_i\mid X_{i-1})
     =  \mathrm{grad} \,F_{K,\mu} (X_{i-1}),
    \qquad i\in\NN,
 \end{align}
 and
  \begin{align*}
   \EE(X_k) = x_0 - \sum_{i=1}^k t_i \EE(\mathrm{grad} \,F_{K,\mu} (X_{i-1})),\qquad k\in\NN.
 \end{align*}
\end{Lem}

\noindent{\bf Proof.}
First note that \ $X_k=x_0-\sum_{i=1}^k t_i Q_i$, \ $k\in\NN$,
 \ where the sequence \ $(Q_i)_{i\in\NN}$ \ is such that the conditional distribution of
 \ $Q_i$ \ with respect to \ $X_{i-1}$ \ is given by
 \begin{align}\label{seged17}
   Q_i =
       \begin{cases}
          \begin{pmatrix}
            1 \\
            1 \\
          \end{pmatrix}
            & \text{with probability \ $\mu(\{(x,y)\in K : X_{i-1}^{(1)}\geq x, \, X_{i-1}^{(2)}\geq y \})$,}\\[1mm]
          \begin{pmatrix}
            1 \\
            -1 \\
          \end{pmatrix}
            & \text{with probability \ $\mu(\{(x,y)\in K : X_{i-1}^{(1)}\geq x, \, X_{i-1}^{(2)}< y \})$,}\\[1mm]
          \begin{pmatrix}
            -1 \\
            1 \\
          \end{pmatrix}
            & \text{with probability \ $\mu(\{(x,y)\in K : X_{i-1}^{(1)}<x, \, X_{i-1}^{(2)}\geq y \})$,}\\[1mm]
         \begin{pmatrix}
            -1 \\
            -1 \\
          \end{pmatrix}
            & \text{with probability \ $\mu(\{(x,y)\in K : X_{i-1}^{(1)}< x, \, X_{i-1}^{(2)}< y \})$.}
      \end{cases}
 \end{align}
Then
 \begin{align*}
   \EE(Q_i\mid X_{i-1})
   &= \begin{pmatrix}
            1 \\
            1 \\
         \end{pmatrix}
         \mu(\{(x,y)\in K : X_{i-1}^{(1)}\geq x, \, X_{i-1}^{(2)}\geq y \})\\
   &\phantom{=\;} + \begin{pmatrix}
              1 \\
              -1 \\
           \end{pmatrix}
          \mu(\{(x,y)\in K : X_{i-1}^{(1)}\geq x, \, X_{i-1}^{(2)}< y \})\\
   &\phantom{=\;}
         + \begin{pmatrix}
            -1 \\
            1 \\
          \end{pmatrix}
     \mu(\{(x,y)\in K : X_{i-1}^{(1)}<x, \, X_{i-1}^{(2)}\geq y \})\\
   &\phantom{=\;}
         + \begin{pmatrix}
            -1 \\
            -1 \\
          \end{pmatrix}
     \mu(\{(x,y)\in K : X_{i-1}^{(1)}< x, \, X_{i-1}^{(2)}< y \})\\
   &=\begin{pmatrix}
       \mu(\{(x,y)\in K : X_{i-1}^{(1)}\geq x\}) - \mu(\{(x,y)\in K : X_{i-1}^{(1)}< x\}) \\
       \mu(\{(x,y)\in K : X_{i-1}^{(2)}\geq y\}) - \mu(\{(x,y)\in K : X_{i-1}^{(2)}< y\}) \\
     \end{pmatrix}
 \end{align*}
 for \ $i\in\NN$.
\ Note that, by Condition (C.3) and Lemma \ref{Lem5}, we also have
  \begin{align*}
   \EE(Q_i\mid X_{i-1})
     =  \begin{pmatrix}
           D_1 F_{K,\mu}(X_{i-1}^{(1)},X_{i-1}^{(2)}) \\
           D_2 F_{K,\mu}(X_{i-1}^{(1)},X_{i-1}^{(2)}) \\
         \end{pmatrix}
     =   \textrm{grad} \,F_{K,\mu} (X_{i-1}),
    \qquad i\in\NN.
 \end{align*}
Hence, by the tower rule, the expectation of \ $X_k$ \ takes the form
 \begin{align*}
   \EE(X_k)  = x_0 - \sum_{i=1}^k t_i \EE(Q_i)
              = x_0 - \sum_{i=1}^k t_i \EE(\EE(Q_i\mid X_{i-1}))
 \end{align*}
 \begin{align*}
             = x_0 - \sum_{i=1}^k t_i \EE(\mathrm{grad} \,F_{K,\mu} (X_{i-1})),
            \qquad k\in\NN.
 \end{align*}
\proofend

\begin{Thm}\label{Thm3}
Let us suppose that Conditions (C.1)--(C.4) hold.
Then the sequence of \ $2$-dimensional random variables defined in \eqref{recursion}
 converges almost surely and in \ $L^q$ \ ($q\in\NN$) \ to the unique global minimizer \ $X^*$ \ of the generalized
 conic function \ $F_{K,\mu}$, \ i.e., \ $\PP(\lim_{n\to\infty} X_n = X^*)=1$ \ and
 \ $\lim_{n\to\infty}\EE\Vert X_n - X^*\Vert^q = 0$.
\end{Thm}

\noindent{\bf Proof.}
First note that under Conditions (C.1)--(C.3) there exists a unique global minimizer \ $\theta^*$ \ of
 \ $F_{K,\mu}$, \ that is, \ $F_{K,\mu}(\theta)>F_{K,\mu}(\theta^*)$ \ for all \ $\theta\ne\theta^*$, $\theta\in\RR^2$,
 \ see, Corollary \ref{Cor1}.
Let us apply Theorem \ref{Thm_Robbins_Monro} with the following choices:
 \begin{itemize}
 \item  $d:=2$, \ $\beta:=0\in\RR^2$, \ and \ $\xi_{n+1}:=Q_{n+1}$, $n\in\ZZ_+$.
  \item $\theta^*\in\RR^2$ \ is such that \ $\mathrm{grad}\ F_{K,\mu}(\theta^*)=0\in\RR^2$.
        \ Note that under the Conditions (C.1)--(C.3), by Corollary \ref{Cor1}, \ $\theta^*$ \ is unique,
        and it is nothing else but the unique global minimizer of \ $F_{K,\mu}$.
 \end{itemize}
In what follows we check that Assumptions (A.1)--(A.3) hold.
Assumption (A.1) holds trivially.
Assumption (A.2) holds with \ $B:=\sqrt{2}$, \ since
 \[
       \left\Vert
        \begin{pmatrix}
            1 \\
            1 \\
        \end{pmatrix}
       \right\Vert
      =
        \left\Vert
        \begin{pmatrix}
            1 \\
            -1 \\
        \end{pmatrix}
       \right\Vert
      =
        \left\Vert
        \begin{pmatrix}
            -1 \\
            1 \\
        \end{pmatrix}
       \right\Vert
      =
       \left\Vert
        \begin{pmatrix}
            -1 \\
            -1 \\
        \end{pmatrix}
       \right\Vert
      =\sqrt{2}.
 \]
Since \ $\EE(Q_i \mid X_0,X_1,\ldots,X_{i-1},Q_1,\ldots,Q_{i-1}) = \EE(Q_i \mid X_{i-1})$,
 \ by \eqref{seged8}, we have \ $M:\RR^2\to\RR^2$, $M(\theta) = \textrm{grad}\ F_{K,\mu}(\theta)$, $\theta\in\RR^2$, \ and,
 by Corollary \ref{Cor1},
 \[
  M(\theta^*) = \mathrm{grad}\ F_{K,\mu}(\theta^*) = 0\in\RR^2.
 \]

Finally, for Assumption (A.3) we have to check that for all \ $\varepsilon\in(0,1)$,
 \[
     \inf_{\varepsilon\leq \Vert \theta - \theta^*\Vert\leq 1/\varepsilon}
            \langle \theta - \theta^*, \mathrm{grad}\ F_{K,\mu}(\theta)\rangle >0.
 \]
Since \ $F_{K,\mu}$ \ is a convex and differentiable function defined on \ $\RR^2$ \ (see, Theorem \ref{Thm2} and
 the proof of Lemma \ref{Lem_conic_uniform}), we have
 \begin{align}\label{seged11}
    \langle \mathrm{grad}\ F_{K,\mu}(\theta), \theta^* - \theta \rangle
      \leq F_{K,\mu}(\theta^*) - F_{K,\mu}(\theta) \leq 0,
      \qquad \forall\;\theta\in\RR^2,
 \end{align}
 where the last inequality follows by that \ $\theta^*$ \ is the global minimizer of \ $F_{K,\mu}$, \ see also Lemma \ref{Lem7}.
Since \ $\theta^*$ \ is strict global minimizer of \ $F_{K,\mu}$, \ i.e.,
 \ $F_{K,\mu}(\theta)>F_{K,\mu}(\theta^*)$ \ for all \ $\theta\ne\theta^*$, $\theta\in\RR^2$ \ (see Corollary \ref{Cor1})
 and \ $\{\theta\in\RR^2 : \varepsilon \leq \Vert \theta - \theta^*\Vert\leq 1/\varepsilon\}$ \ is a compact set,
 by \eqref{seged11}, we get Assumption (A.3) holds in our case.
\proofend

\begin{Ex}
Let \ $K$ \ be the square with vertexes \ $(0,0), (0,1), (1,0), (1,1)$ \  as in part (i) of Example \ref{example1}.
Let us assume that \ $\mu$ \ is the probability measure on \ $K$ \ with Radon-Nikodym derivative with respect to the
 Lebesgue measure given by
 \[
   h_\mu(x,y)=\begin{cases}
                   1 & \text{if \ $(x,y)\in K$,}\\
                   0 & \text{if \ $(x,y)\not\in K$.}
                 \end{cases}
 \]
Further, let \ $x_0:=(0, 0)^\top$ \ and \ $t_k:=\frac{1}{k}$, $k\in\NN$.
\ Then
 \begin{align*}
   X_0=\begin{pmatrix}
            0 \\
            0 \\
       \end{pmatrix},
   \qquad
   X_k=-\sum_{i=1}^k t_i Q_i = -\sum_{i=1}^k \frac{1}{i} Q_i ,\quad k\in\NN,
 \end{align*}
 where the sequence \ $(Q_i)_{i\in\NN}$ \ is such that the conditional distribution of \ $Q_i$ \
 with respect to \ $X_{i-1}$ \ is given by \eqref{seged17}.
By Theorem \ref{Thm3} and part (i) of Example \ref{example1}, we have \ $\PP(\lim_{k\to\infty} X_k = X^*)=1$ \ and
 \ $\lim_{k\to\infty}\EE\Vert X_k - X^*\Vert^q = 0$ \ for all \ $q\in\NN$, \ where \ $X^*=(1/2 , 1/2)^\top$.
\ Note also that if \ $X_{i-1}\in K$, \ then the conditional distribution of \ $Q_i$ \ with respect to \ $X_{i-1}$ \ takes the form
 \begin{align*}
   Q_i =
       \begin{cases}
          \begin{pmatrix}
            1 \\
            1 \\
          \end{pmatrix}
            & \text{with probability \ $X_{i-1}^{(1)}X_{i-1}^{(2)}$,}\\
          \begin{pmatrix}
            1 \\
            -1 \\
          \end{pmatrix}
            & \text{with probability \ $X_{i-1}^{(1)}(1-X_{i-1}^{(2)})$,}\\
          \begin{pmatrix}
            -1 \\
            1 \\
          \end{pmatrix}
            & \text{with probability \ $(1-X_{i-1}^{(1)})X_{i-1}^{(2)}$,}\\
         \begin{pmatrix}
            -1 \\
            -1 \\
          \end{pmatrix}
            & \text{with probability \ $(1-X_{i-1}^{(1)})(1-X_{i-1}^{(2)})$.}
      \end{cases}
 \end{align*}
Finally, we remark that \ $X_1=(1,1)^\top$ \ and \ $X_2=(1/2,1/2)^\top$.
\end{Ex}

\subsection{Almost sure and \ $L^q$-convergence of \ $(F_{K,\mu}(X_k))_{k\in\ZZ_+}$}

First we recall an equivalent reformulation of \ $L^q$-convergence, where \ $q\in\NN$, \ see, e.g.,
 Chow and Teicher \cite[Theorem 4.2.3]{ChoTei}.

\begin{Lem}\label{Lem_Lq}
Let \ $d,q\in\NN$, \ $\xi:\Omega\to\RR^d$ \ and \ $\xi_n:\Omega\to\RR^d$, $n\in\NN$, \ be \ $\RR^d$-valued random variables such
 that \ $\EE(\Vert \xi\Vert^q)<\infty$ \ and \ $\EE(\Vert \xi_n\Vert^q)<\infty$, $n\in\NN$.
\ Then \ $\xi_n$ \ converges to \ $\xi$ \ in \ $L^q$ \ as \ $n\to\infty$ \ (i.e., \ $\lim_{n\to\infty} \EE(\Vert \xi_n-\xi\Vert^q)=0$)
 \  if and only if \ $\xi_n$ \ converges in probability to \ $\xi$ \ as \ $n\to\infty$ \ and the set of random variables
 \ $\{ \Vert \xi_n\Vert^q : n\in\NN\}$ \ is uniformly integrable, i.e.,
 \[
    \lim_{m\to\infty} \sup_{n\in\NN} \EE\left( \Vert \xi_n\Vert^q \mathbf 1_{\{ \Vert \xi_n\Vert^q >m\}}\right)=0.
 \]
\end{Lem}

\begin{Thm}\label{Thm4}
Let us suppose that Conditions (C.1)--(C.4) hold.
Then the sequence of one-dimensional random variables \ $(F_{K,\mu}(X_k))_{k\in\NN}$ \ converges almost surely
 and in \ $L^q$ \ ($q\in\NN$) \ to \ $F_{K,\mu}(X^*)$ \ as \ $k\to\infty$, \ where \ $X^*$ \ denotes the unique
 global minimizer of \ $F_{K,\mu}$.
\end{Thm}

\noindent{\bf Proof.}
By Theorem \ref{Thm3}, \ $\PP(\lim_{k\to\infty}X_k = X^*)=1$, \ and hence to prove that
 \ $\PP(\lim_{k\to\infty} F_{K,\mu}(X_k) = F_{K,\mu}(X^*))=1$, \ it is enough to check that \ $F_{K,\mu}$ \ is continuous.
This follows by that \ $F_{K,\mu}$ \ is a convex function defined on \ $\RR^2$ \ (see Theorem \ref{Thm2}).
We give an alternative argument, too.
Let \ $(x_n,y_n)^\top\in\RR^2$, $n\in\NN$, \ be such that \ $\lim_{n\to\infty} (x_n,y_n) = (x,y)$, \ where
 \ $(x,y)^\top\in\RR^2$.
\ Then for all \ $(\alpha,\beta)^\top\in\RR^2$, \ $\lim_{n\to\infty}d_1((x_n,y_n),(\alpha,\beta)) = d_1((x,y),(\alpha,\beta))$,
 \ and, using that \ $K$ \ is bounded,
 \[
    \sup_{n\in\NN} \sup_{(\alpha,\beta)\in K} d_1((x_n,y_n),(\alpha,\beta)) <\infty.
 \]
By Lebesgue dominated convergence theorem (which can be used since \ $\mu(K)<\infty$)
 \begin{align*}
  \lim_{n\to\infty} F_{K,\mu}(x_n,y_n)
    &= \int_K \lim_{n\to\infty} d_1((x_n,y_n),(\alpha,\beta)) \,\mu(\dd\alpha,\dd\beta)\\
    & = \int_K d_1((x,y),(\alpha,\beta)) \,\mu(\dd\alpha,\dd\beta)
      = F_{K,\mu}(x,y),
 \end{align*}
 yielding that \ $F_{K,\mu}$ \ is continuous.

\begin{figure}[h!]
\centering
\includegraphics[height=3cm]{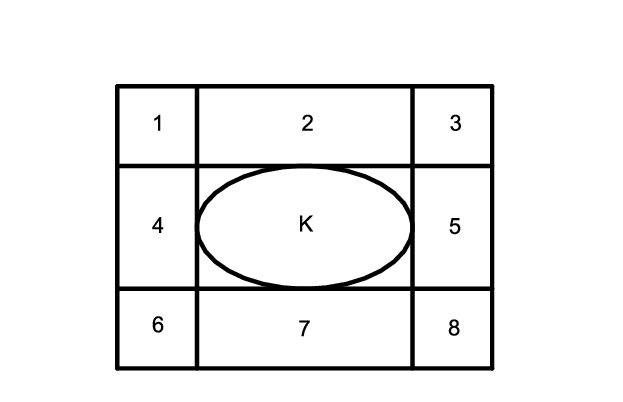}
\caption{The eight cases.}
\label{figure2}
\end{figure}

Further, using Lemma \ref{Lem_Lq} and that almost sure convergence yields convergence in probability, in order to prove
 \ $L^q$-convergence of  \ $(F_{K,\mu}(X_k))_{k\in\NN}$, \ it is enough (and actually necessary) to check that
 \begin{align}\label{seged18}
    \lim_{m\to\infty} \sup_{k\in\NN} \EE\left( \Vert X_k\Vert^q \mathbf 1_{\{ \Vert X_k\Vert^q >m\}}\right)=0.
 \end{align}
We show that the sequence \ $(\Vert X_k\Vert^q)_{k\in\NN}$ \ is bounded, and then \eqref{seged18} readily follows.
Let \ $D:=\sup_{k\in\NN}\{t_k\}=t_1>0$ \ (indeed, \ $(t_k)_{k\in\NN}$ \ is a decreasing sequence of positive numbers).
Let us consider the rectangle \ $R$ \ with vertexes
 \begin{align*}
  &\left( \inf\{x : (x,y)\in K\} - D\sqrt{2}, \; \inf\{y : (x,y)\in K\} - D\sqrt{2} \right),\\
  &\left( \inf\{x : (x,y)\in K\} - D\sqrt{2}, \; \sup\{y : (x,y)\in K\} + D\sqrt{2} \right),\\
  &\left( \sup\{x : (x,y)\in K\} + D\sqrt{2}, \; \inf\{y : (x,y)\in K\} - D\sqrt{2} \right),\\
  &\left( \sup\{x : (x,y)\in K\} + D\sqrt{2}, \; \sup\{y : (x,y)\in K\} + D\sqrt{2} \right).
 \end{align*}
Since \ $\Vert Q_k\Vert = \sqrt{2}$, $k\in\NN$, \ if \ $X_n\in K$ \ with some \ $n\in\NN$, \ then \ $X_{n+1}\in R$,
 \ i.e., the recursion \eqref{recursion} cannot leave the rectangle \ $R$ \ starting from \ $K$ \ by one step.
Next we check that if \ $X_n\in R$ \ with some \ $n\in\NN$, \ then \ $X_{n+1}\in R$, \ which yields that
 the recursion \eqref{recursion} cannot leave the rectangle \ $R$.
\ We distinguish eight cases according to the Figure \ref{figure2}.
If \ $X_n$ \ is in the rectangle numbered 1, then \ $Q_{n+1} = (-1,1)^\top$ \ and hence, by the choice of \ $D$,
 \[
    X_{n+1} = X_n + t_{n+1}\begin{pmatrix}
                             1 \\
                             -1 \\
                           \end{pmatrix}
                           \in R.
 \]
If \ $X_n$ \ is in the rectangle numbered 2, then \ $Q_{n+1} = (1,1)^\top$ \ or \ $Q_{n+1} = (-1,1)^\top$ \ according to the cases
 \ $X_n^{(1)}\geq P_{n+1}^{(1)}$ \ and \ $X_n^{(1)}< P_{n+1}^{(1)}$, \ and hence
 \[
    X_{n+1} = X_n + t_{n+1}\begin{pmatrix}
                             -1 \\
                             -1 \\
                           \end{pmatrix}
                           \in R
   \qquad \text{or}\qquad
    X_{n+1} = X_n + t_{n+1}\begin{pmatrix}
                             1 \\
                             -1 \\
                           \end{pmatrix}
                           \in R.
  \]
If \ $X_n$ \ is in the rectangle numbered 3, then \ $Q_{n+1} = (1,1)^\top$ \ and hence
 \[
    X_{n+1} = X_n + t_{n+1}\begin{pmatrix}
                             -1 \\
                             -1 \\
                           \end{pmatrix}
                           \in R.
 \]
The other cases can be handled similarly.
\proofend

\section*{Acknowledgements}
We are grateful for the referee for his/her several valuable comments
 that have led to an improvement of the manuscript.

\end{document}